\pdfoutput=1
\RequirePackage{ifpdf}
\ifpdf 
\documentclass[pdftex]{sigma}
\else
\documentclass{sigma}
\fi

\def\sw#1{{\sb{(#1)}}}
\newcommand{\PW}{{\mathcal{P}_{C({\rm U}(1))}(C(S^3))}}
\newcommand{\PWsq}{{\PW\square_{\mathcal{O}({\rm U}(1))}\mathcal{O}({\rm SU}_q(2))}}

\begin{document}

\allowdisplaybreaks

\renewcommand{\thefootnote}{$\star$}

\renewcommand{\PaperNumber}{031}

\FirstPageHeading

\ShortArticleName{Nontrivial Deformation of a~Trivial Bundle}

\ArticleName{Nontrivial Deformation of a~Trivial Bundle\footnote{This paper is a~contribution to the Special Issue on
Noncommutative Geometry and Quantum Groups in honor of Marc A.~Rief\/fel.
The full collection is available at
\href{http://www.emis.de/journals/SIGMA/Rieffel.html}{http://www.emis.de/journals/SIGMA/Rieffel.html}}}

\Author{Piotr~M.~HAJAC~$^{\dag\ddag}$ and Bartosz~ZIELI\'NSKI~$^\S$}

\AuthorNameForHeading{P.M.~Hajac and B.~Zieli\'nski}

\Address{$^\dag$~Instytut Matematyczny, Polska Akademia Nauk, ul.~\'Sniadeckich 8, 00-956 Warszawa, Poland}
\EmailD{\href{mailto:pmh@impan.pl}{pmh@impan.pl}} \URLaddressD{\url{http://www.impan.pl/~pmh/}}
\Address{$^\ddag$~Katedra Metod Matematycznych Fizyki, Uniwersytet Warszawski,\\
\hphantom{$^\ddag$}~ul.~Ho\.za~74, 00-682 Warszawa, Poland}

\Address{$^\ddag$~Department of Computer Science, Faculty of Physics and Applied Informatics,\\
\hphantom{$^\S$}~University of \L{}\'od\'z, Pomorska 149/153 90-236 \L{}\'od\'z, Poland}
\EmailD{\href{mailto:bzielinski@uni.lodz.pl}{bzielinski@uni.lodz.pl}}

\ArticleDates{Received October 29, 2013, in f\/inal form March 03, 2014; Published online March 27, 2014}

\Abstract{The ${\rm SU}(2)$-prolongation of the Hopf f\/ibration $S^3\to S^2$ is a~trivializable principal ${\rm SU}(2)$-bundle.
We present a~noncommutative deformation of this bundle to a~quantum principal ${\rm SU}_q(2)$-bundle that is \emph{not}
trivializable.
On the other hand, we show that the ${\rm SU}_q(2)$-bundle is \emph{piecewise} trivializable with respect to the closed
covering of $S^2$ by two hemispheres intersecting at the equator.}


\Classification{58B32}

\rightline{\it Dedicated to Marc A.~Rieffel on the occasion of his 75th birthday}

\renewcommand{\thefootnote}{\arabic{footnote}} 
\setcounter{footnote}{0}

\section{Introduction and preliminaries}

The goal of this paper is to show how a~noncommutative deformation can turn a~trivializable principal bundle into
a~nontrivializable quantum principal bundle.
This is a~peculiar phenomenon because noncommutative deformations usually preserve basic topological features of
deformed objects, e.g.~$K$-groups.

On the other hand, this paper exemplif\/ies the general theory of piecewise trivial principal comodule algebras developed
in~\cite{hkmz11,hrz}.
Therefore we follow the notation, conventions and general setup employed therein.
To make our exposition self-contained and easy to read, we often recall basic concepts and def\/initions.

Let $\pi:X\rightarrow M$ be a~principal $G$-bundle over $M$, and $G'$ be a~subgroup of $G$.
A $G'$-{\em reduction} of $X\rightarrow M$ is a~subbundle $X'\subseteq X$ over $M$ that is a~principal $G'$-bundle over
$M$ via the restriction of the $G$-action on $X$.
Many important structures on manifolds can be formulated as reductions of their frame bundles.
For instance, an orientation, a~volume form and a~metric on a~manifold $M$ correspond to reductions of the frame bundle
$FM$ to a~${\rm GL}_+(n,{\mathbb R})$, ${\rm SL}(n,{\mathbb R})$ and ${\rm O}(n,{\mathbb R})$-bundle, respectively.
See~\cite{kn63} for more details.

An operation inverse to a~reduction of a~principal bundle is a~prolongation of a~principal bundle.
Let $\pi:X\rightarrow M$ be a~principal $G'$-bundle over $M$, and let $G'$ be a~subgroup of $G$.
Def\/ine $X\times_{G'} G:=(X\times G)/{\sim}$, where $(x,g)\sim (xh,h^{-1}g)$, for all $x\in X$, $g\in G$ and $h\in G'$.
Then
\begin{gather*}
\hat\pi: \ X\times_{G'} G \longrightarrow   M,
\qquad
[x,g] \longmapsto   \pi(x),
\end{gather*}
is a~$G$-bundle called the $G$-\emph{prolongation} of $X$, with the $G$-action given by $[x,g]h:=[x,gh]$.
The bundle $X\rightarrow M$ is a~$G'$-reduction of $X\times_{G'} G\rightarrow M$.

An interesting special case is when $X=G$ and $M=G/G'$, that is the homogenous bundle case.
It is easy to see that $G\times_{G'}G\rightarrow G/G'$ is always a~trivializable bundle.
Indeed, the following $G$-equivariant bundle maps provide an explicit isomorphism and its inverse:
\begin{gather*}
f: \ G\times_{G'}G \longrightarrow   G/G'\times G,
\qquad
[g_1,g_2] \longmapsto   ([g_1],g_1g_2),
\nonumber
\\
f^{-1}: \ G/G'\times G \longrightarrow   G\times_{G'}G,
\qquad
([g],h) \longmapsto   [g,g^{-1}h].
\end{gather*}

A quantum-group version of the trivializability of $G\times_{G'}G\rightarrow G/G'$ can be easily checked mimicking the
classical argument.
In particular, the ${\rm SU}_q(2)$-prolongation
\begin{gather*}
{\rm SU}_q(2)\underset{{\rm U}(1)}{\times}{\rm SU}_q(2) \longrightarrow   S^2_q
\end{gather*}
of the standard quantum Hopf f\/ibration is trivializable~\cite[p.~1104]{bz12}.
However, as the main result of this paper, we show that the ${\rm SU}_q(2)$-prolongation
\begin{gather*}
{\rm SU}(2)\underset{{\rm U}(1)}{\times}{\rm SU}_q(2) \longrightarrow   S^2
\end{gather*}
of the classical Hopf f\/ibration is \emph{not} trivializable.

\subsection{Notation}

\looseness=-1
We work over the f\/ield~$\mathbb{C}$ of complex numbers.
The unadorned tensor product stands for the tensor product over this f\/ield.
The comultiplication, counit and the antipode of a~Hopf algebra~$H$ are denoted by~$\Delta$, $\varepsilon$ and $S$,
respectively.
Our standing assumption is that $S$ is invertible.
A~right $H$-comodule algebra $P$ is a~unital associative algebra equipped with an $H$-coaction $ \Delta_P: P \rightarrow
P \otimes H$ that is an algebra homomorphism.
For a~comodule algebra $P$, we call
\begin{gather*}
P^{\operatorname{co} H}:=\left\{p \in P\,|\,\Delta_P(p)=p \otimes 1\right\}
\end{gather*}
the subalgebra of coaction-invariant elements in~$P$.
A left coaction on $V$ is denoted by ${}_V\Delta$.
For comultiplications and coactions, we often employ the Heynemann--Sweedler notation with the summation symbol
suppressed:
\begin{gather*}
\Delta(h)=:h_{(1)}\otimes h_{(2)},
\qquad
\Delta_P(p)=:p_{(0)}\otimes p_{(1)},
\qquad
{}_V\Delta(v)=:v_{(-1)}\otimes v_{(0)}.
\end{gather*}

\subsection{Reductions and prolongations of principal comodule algebras}

\begin{definition}[\cite{bh04}]
 Let $H$ be a~Hopf algebra, $P$ be a~right $H$-comodule algebra and let $B:=P^{\operatorname{co} H}$ be the
coaction-invariant subalgebra.
The comodule algebra $P$ is called {\em principal}~if\/f:
\begin{enumerate}\itemsep=0pt
\item[1)]
$P{\otimes}_B P\ni p \otimes q \mapsto\operatorname{can}(p \otimes q):= pq_{(0)} \otimes q_{(1)}\in P \otimes H$ is
bijective,
\item[2)]
there exists a left $B$-linear right $H$-colinear splitting of the multiplication map $B\otimes P\to P$,
\item[3)]
the antipode of $H$ is bijective.
\end{enumerate}
\end{definition}

Here (1) is the Hopf--Galois (freeness) condition, (2) means equivariant projectivity of $P$, and (3) ensures
a~left-right symmetry of the def\/inition (everything can be re-written for left comodule algebras).

A particular class of principal comodule algebras is distinguished by the existence of a~clea\-ving map.
A cleaving map is def\/ined as a~unital right $H$-colinear convolution-invertible map \mbox{$j:H\rightarrow P$}.
Comodule algebras admitting a~cleaving map are called {\em cleft}.
One can show that a~cleaving map is automatically injective.
However, in general, they are not algebra homomorphisms.

If $j:H\rightarrow P$ is a~right $H$-colinear algebra homomorphism, then it is automatically convo\-lu\-tion-invertible
and unital.
A cleft comodule algebra admitting a~cleaving map that is an algebra homomorphism is called a~{\em smash product}.
All commutative smash products reduce to the tensor algebra $P^{\operatorname{co} H}\otimes H$, so that smash products
play the role of trivial bundles.
Here a~cleaving map is simply given by $j(h):=1\otimes h$.
A cleaving map def\/ines a~left action of $H$ on $P^{\operatorname{co} H}$ making it a~left $H$-module algebra:
$h\triangleright p:=j(h_{(1)})pj^{-1}(h_{(2)})$.
Conversely, if $B$ is a~left $H$-module algebra, one can construct a~smash product $B\rtimes H$ by equipping the vector
space $B\otimes H$ with the multiplication
\begin{gather*}
(a\otimes h)(b\otimes k):=a\,(h_{(1)}\triangleright b)\otimes h_{(2)}\,k,
\qquad
a,b\in B,
\qquad
h,k\in H,
\end{gather*}
and coaction $\Delta_{B\rtimes H}:={\rm id}\otimes\Delta$.
Then again a~cleaving map is simply given by $j(h):=1\otimes h$.

\begin{definition}[\cite{g-r99,qsng,s-p99}]
Let $P$ be a~principal $H$-comodule algebra and $J$ be a~Hopf ideal of $H$ such that $H$ is a~principal left
$H/J$-comodule algebra.
We say that an ideal $I$ of $P$ is a~$J$-\emph{reduction} of $P$ if and only if the following conditions are satisf\/ied:
\begin{enumerate}\itemsep=0pt
\item[1)]
$I$ is an $H/J$-subcomodule of $P$,
\item[2)]
$P/I$ with the induced coaction is a~principal $H/J$-comodule algebra,
\item[3)]
$(P/I)^{\operatorname{co} H/J}=P^{\operatorname{co} H}$.
\end{enumerate}
\end{definition}

Loosely speaking, $J$ plays the role of the ideal of functions vanishing on a~subgroup and~$I$ the ideal of
functions vanishing on a~subbundle.
Thus~$H/J$ works as the algebra of the reducing subgroup, and $P/I$ as the algebra of the reduced bundle.
The coaction-invariant subalgebra~$P^{\operatorname{co} H}$ remains intact~-- the base space of a~subbundle coincides
with the base space of the bundle.

If $M$ is a~right comodule over a~coalgebra $C$ and $N$ is a~left $C$-comodule, then we def\/ine their {\em cotensor
product} as
\begin{gather*}
M\underset{C}{\Box}N:=\{t\in M\otimes N\;|\;(\Delta_M\otimes{\rm id})(t)= ({\rm id}\otimes{}_N\Delta)(t)\}.
\end{gather*}
In particular, for a~principal $H'$-comodule algebra $P$ and a~Hopf algebra epimorphism $H\stackrel{\pi}{\rightarrow}
H'$ making $H$ a~left $H'$-comodule in the obvious way, one proves that the cotensor product $P\Box_{H'}H$ is
a~principal $H$-comodule algebra with the  $H$-coaction def\/ined by ${\rm id}\otimes\Delta$.
We call the principal comodule algebra $P\Box_{H'}H$ the \emph{$H$-prolongation of $P$}.

\subsection{Piecewise triviality}

\begin{definition}[\protect{cf.\ \cite[Def\/inition~3.6]{hkmz11}}] A family of surjective algebra homomorphisms\linebreak $\{\pi_i:P\rightarrow
P_i\}_{i\in\{1,\ldots,N\}}$, $N\in\mathbb{N}\setminus\{0\}$, is called a~{\em covering} if\/f
\begin{enumerate}\itemsep=0pt
\item[1)]
$\bigcap_{i\in\{1,\ldots,N\}}\operatorname{Ker}\pi_i=0$,
\item[2)]
The family of ideals $(\operatorname{Ker}\pi_i)_{i\in\{1,\ldots,N\}}$ generates a~distributive lattice with $+$ and
$\cap$ as meet and join respectively.
\end{enumerate}
\end{definition}

We recall now (cf.~\cite[Def\/inition~3.8]{hkmz11}) a~quantum version of the notion of piecewise triviality of principal
bundles (like local triviality, but with respect to closed subsets).
\begin{definition}
An $H$-comodule algebra $P$ is called \emph{piecewise trivial} if\/f there exists a~family $\{\pi_i:P\rightarrow
P_i\}_{i\in\{1,\ldots,N\}}$, $N\in\mathbb{N}\setminus\{0\}$, of surjective $H$-colinear maps such that:
\begin{enumerate}\itemsep=0pt
\item[1)]
the restrictions $\pi_i|_{P^{\operatorname{co} H}}:P^{\operatorname{co} H}\rightarrow P_i^{\operatorname{co} H}$ form
a~covering,
\item[2)]
the $P_i$'s are smash products ($P_i\cong P_i^{\operatorname{co} H}\rtimes H$ as $H$ comodule algebras).
\end{enumerate}
\end{definition}

Assume also that the antipode of $H$ is bijective.
Then, as smash products are principal, it follows from~\cite[Theorem~3.3]{hkmz11} that piecewise trivial comodule
algebras are automatically principal.
To emphasize this fact and stay in touch with the classical terminology, we frequently use the phrase ``piecewise
trivial principal comodule algebra''.
Note also that the consequence of principality of $P$ is that $\{\pi_i:P\rightarrow P_i\}_{i\in\{1,\ldots,N\}}$ is
a~covering of $P$ (see~\cite{hrz}).

\begin{definition}[\cite{hrz}]\sloppy
Let $\{\pi_i:P\to P_i\}_{i\in\{1,\ldots,N\}}$, $N\in\mathbb{N}\setminus\{0\}$, be a~covering by
right  $H$-co\-li\-near maps of a~principal right $H$-comodule algebra $P$ such that the restrictions
\mbox{$\pi_i|_{P^{\operatorname{co} H}}:P^{\operatorname{co} H}\rightarrow P_i^{\operatorname{co} H}$} also form a~covering.
A {\em piecewise trivialization} of $P$ with respect to the covering \mbox{$\{\pi_i:P\to P_i\}_{i\in\{1,\ldots,N\}}$} is
a~family $\{j_i:H\to P_i\}_{i\in\{1,\ldots,N\}}$ of right $H$-colinear algebra homomorphisms (cleaving maps).
\end{definition}

It is clear that a~principal comodule algebra is piecewise trivial if and only if it admits a~piecewise
trivialization (see the preceding section).

\subsection{The Peter--Weyl comodule algebra}

The Peter--Weyl comodule algebra (see~\cite{bdh} and references therein) extends the notion of re\-gular functions in the
$C^*$-algebra of a~compact quantum group (linear combinations of matrix coef\/f\/icients of the f\/inite-dimensional
corepresentations) to unital $C^*$-algebras equipped with a~compact quantum group action.

\begin{definition}[cf.~\cite{p-p95}] For a~unital $C^*$-algebra $A$ and a~compact quantum group $(H,\Delta)$, we say that an injective
unital $*$-homomorphism $\delta:A\rightarrow A\otimes_{\min}H$ is a~\emph{coaction} if and only if
\begin{enumerate}\itemsep=0pt
\item[1)]
$(\delta\otimes\mathrm{id})\circ\delta= (\mathrm{id}\otimes\Delta)\circ\delta$ (coassociativity),
\item[2)]
$\{\delta(a)(1\otimes h)\, |\, a\in A,\,h\in H\}^{\mathrm{cls}}= A\underset{\min}{\otimes}H$ (counitality).
\end{enumerate}
  Here $\otimes_{\min}$ denotes the spatial tensor product of $C^*$-algebras and
$\{\cdot\}^{\mathrm{cls}}$ stands for the closed linear span of a~subset of a~Banach space.
We say that a~compact quantum group \emph{acts} on a~unital $C^*$-algebra if there is a~coaction in the aforementioned
sense.
\end{definition}

Next, we denote by ${\mathcal{O}} (H)$ the dense Hopf $*$-subalgebra of $H$ spanned by the matrix coef\/f\/i\-cients of
f\/inite-dimensional corepresentations.
We def\/ine the \emph{Peter--Weyl subalgebra} of $A$~\cite{bdh} as
\begin{gather*}
{\mathcal{P}}_H(A):=\{ a\in A\,| \,\delta(a)\in A\otimes{\mathcal{O}} (H) \}.
\end{gather*}
One shows that it is an ${\mathcal{O}} (H)$-comodule algebra which is a~dense $*$-subalgebra of~$A$~\cite{p-p95,s-pm11}.

The Peter--Weyl comodule algebra of functions on a~compact Hausdorf\/f space with an action of a~compact group is principal
if and only if the action is free~\cite{bdh, bh}.
In other words, the Galois condition of Hopf--Galois theory holds if and only if we have a~compact principal bundle.

\section[The ${\rm SU}_q(2)$-prolongation of the classical Hopf f\/ibration]{The $\boldsymbol{{\rm SU}_q(2)}$-prolongation of the classical Hopf f\/ibration}

To f\/ix the notation, let us recall def\/initions of the Hopf algebras $\mathcal{O}({\rm U}(1))$ and $\mathcal{O}({\rm SU}_q(2))$, and
the Peter--Weyl comodule algebra $\PW$ of functions on the classical sphere~$S^3$.
For details on the latter algebra we refer the reader to~\cite{bhms07}.

Recall that the $^*$-algebra $\mathcal{O}({\rm U}(1))$ of polynomial functions on ${\rm U}(1)$ is generated by the unitary element
$u:{\rm U}(1)\ni x\mapsto x\in {\mathbb C}$, and can be equivalently def\/ined as the algebra of Laurent polynomials in $u$ subject
to the relation $u^{-1}=u^*$.
The Hopf algebra structure is given by $\Delta(u):=u\otimes u$, $\varepsilon(u):=1$ and $S(u):=u^{-1}$.

The algebra of polynomial functions on ${\rm SU}_q(2)$~\cite{w-sl87} is generated as a~${}^*$-algebra by $\alpha$ and $\gamma$
satisfying relations
\begin{gather}
\label{ac}
\alpha\gamma=q\gamma\alpha,
\qquad
\alpha\gamma^*=q\gamma^*\alpha,
\qquad
\gamma\gamma^*=\gamma^*\gamma,
\qquad
\alpha^*\alpha+\gamma^*\gamma=1,
\qquad
\alpha\alpha^*+q^2\gamma\gamma^*=1,\!\!\!
\end{gather}
where $0<q\leq 1$.
The Hopf algebra structure comes from the matrix
\begin{gather*}
U:=
\begin{pmatrix}
\alpha & -q\gamma^*
\\
\gamma &\alpha^*
\end{pmatrix},
\qquad
\text{i.e.}
\qquad
\Delta(U_{ij}):=\sum\limits_k U_{ik}\otimes U_{kj},
\quad
S(U_{ij}):=U_{ji}^*,
\quad
\varepsilon(U_{ij}):=\delta_{ij}.
\end{gather*}
The Hopf ${}^*$-algebra epimorphism
\begin{gather}
\label{mudef}
\pi: \ \mathcal{O}({\rm SU}_q(2)) \longrightarrow  \mathcal{O}({\rm U}(1)),
\qquad
\pi(\alpha):=u,
\qquad
\pi(\gamma):=0,
\end{gather}
makes $\mathcal{O}({\rm SU}_q(2))$ into a~left and right $\mathcal{O}({\rm U}(1))$-comodule algebra via the left and right coactions
$(\pi\otimes{\rm id})\circ\Delta$ and $({\rm id}\otimes\pi)\circ\Delta$ respectively.
For $q=1$ the Hopf algebra $\mathcal{O}({\rm SU}_q(2))$ is commutative, and we denote its generators by $a$ and $c$ rather
then $\alpha$ and $\gamma$.

The Peter--Weyl comodule algebra $\PW$ is the subalgebra of $C({\rm SU}(2))$ that is the \emph{algebraic} direct sum of the
modules of continuous sections of the complex line bundles $L_n$, $n\in\mathbb{Z}$, associated to the Hopf f\/ibration:
\begin{gather*}
\PW=\bigoplus_{n\in{\mathbb Z}}\Gamma(L_n).
\end{gather*}
We have the following proper inclusions of function algebras:
\begin{gather*}
\mathcal{O}({\rm SU}(2))\subsetneq \PW\subsetneq C(S^3).
\end{gather*}

Next, recall that $a,c:S^3\rightarrow {\mathbb C}$ are coordinate functions on $S^3$ satisfying $|a|^2+|c|^2=1$.
The diagonal action of ${\rm U}(1)$ on $S^3$ yielding the Hopf f\/ibration dualizes to the ${\mathcal{O}}({\rm U}(1))$-comodule
algebra structure on $\PW$ given by $a\mapsto a\otimes u$, $c\mapsto c\otimes u$.

Now we will describe the piecewise trivial structure of $\PW$.
For brevity, we def\/ine
\begin{gather*}
\omega:=\sqrt{\frac{2}{1+||a|^2-|c|^2|}}.
\end{gather*}
Note that $\omega$ is an element of the coaction-invariant subalgebra
$\mathcal{P}_{C(\mathrm{U}(1))}(C(S^3))^{\operatorname{co} \mathcal{O}(\mathrm{U}(1))}
=\Gamma(L_0)$, which we identify with $C(S^2)$.
Let us also def\/ine the following ideals $I_a,I_c\subseteq C(S^2)$:
\begin{gather*}
I_a:=\big\{f\in C(S^2)\,\big|\,f(x)=0\text{\;for all\;}x\in S^2\;\text{such that\;}|a|^2(x)\leq 1/2\big\},
\nonumber
\\
I_c:=\big\{f\in C(S^2)\,\big|\, f(x)=0\text{\;for all\;}x\in S^2\;\text{such that\;}|a|^2(x)\geq 1/2\big\}.
\end{gather*}
It is well known (cf.~\cite{bhms07}) that the canonical surjections $C(S^2)\to C(S^2)/I_i\cong C(D)$, $i=a,c$, where $D$ is the unit disk, form a covering,
and that $(1-\omega^2|a|^2)\in I_a$, $(1-\omega^2|c|^2)\in I_c$.
We also know \cite[equation~(3.4.57)]{bhms07} that
\begin{gather*}
\big(1-\omega^2|a|^2\big)\big(1-\omega^2|c|^2\big)=0.
\end{gather*}

The covering of $\PW$ can now be given by the canonical surjections in terms of $I_a$ and $I_c$ (cf.~\cite{bhms07}):
\begin{gather*}
\begin{split}
& \pi_a: \ \PW \longrightarrow  \PW/(I_a\PW),
\\
& \pi_c: \ \PW \longrightarrow  \PW/(I_c\PW).
\end{split}
\end{gather*}
Indeed, since $\mathcal{P}_{C(\mathrm{U}(1))}(C(S^3))$ is a principal $\mathcal{O}(\mathrm{U}(1))$-comodule algebra with the coaction-inva\-riant subalgebra $C(S^2)$, it follows from \cite[Proposition~3.4]{hkmz11}  that  the maps $\pi_i$ form a covering.

A trivialization associated with the above covering is given by the following cleaving maps, which are clearly algebra
homomorphisms:
\begin{gather*}
j_a:\  \mathcal{O}({\rm U}(1)) \longrightarrow  \PW/(I_a\PW),
\qquad
u^n \longmapsto   \pi_a(\omega a)^n,
\\
j_c:\  \mathcal{O}({\rm U}(1)) \longrightarrow \PW/(I_c\PW),
\qquad
u^n \longmapsto  \pi_c(\omega c)^n.
\end{gather*}
One can argue (cf.~\cite{bhms07}) that
\begin{gather}
f_a:\  \PW/(I_a\PW)\stackrel{\cong}{\longrightarrow} C(D)\otimes\mathcal{O}({\rm U}(1)),
\nonumber
\\
f_c:\ \PW/(I_c\PW)\stackrel{\cong}{\longrightarrow} C(D)\otimes\mathcal{O}({\rm U}(1)),
\label{sphids}
\\
f_i: \ x  \longmapsto   \pi_i(x\sw{0})j_i(Sx\sw{1})\otimes x\sw{2},
\qquad
i=a,c.
\nonumber
\end{gather}
To see this, f\/irst note that $\pi_a(C(S^2))\cong C(D)\cong \pi_c(C(S^2))$.
Then, for any $n\in{\mathbb Z}$,
\begin{gather*}
f_i(\pi_i(\Gamma(L_n)))=\pi_i (C(S^2) )\otimes u^n=C(D)\otimes u^n,
\qquad
i=a,c,
\end{gather*}
whence $f_i(\PW)=C(D)\otimes \mathcal{O}({\rm U}(1))$.
Indeed, $f_i(\pi_i(\Gamma(L_n)))\subseteq\pi_i(C(S^2))\otimes u^n$.
On the other hand, consider an arbitrary element $y\in C(D)$.
Then there exist elements $y_a,y_c\in C(S^2)$ such that $y=\pi_a(y_a)=\pi_c(y_c)$.
Hence
\begin{gather*}
y\otimes u^n=f_z(\pi_z(y_zz^n\omega^n)),
\qquad
z=a,c.
\end{gather*}
Here we adopt the convention that $z^{-|n|}:=(z^*)^{|n|}$. Summarizing, $\mathcal{P}_{C(\mathrm{U}(1))}(C(S^3))$ is a piecewise trivial principal comodule algebra~\cite{bhms07}.

Since $\mathcal{O}({\rm SU}_q(2))$ is a~left principal $\mathcal{O}({\rm U}(1))$-comodule algebra, by~\cite[Lemma~1.13]{hrz} the
cotensor product $\PWsq$ is a~piecewise trivial principal comodule algebra.
Explicitly, the covering and trivializations inherited from $\PW$ make it piecewise trivial via the formulas:
\begin{gather*}
\hat{\pi}_i:=\pi_i\otimes{\rm id},
\qquad
\hat{j}_i:=(j_i\circ\pi\otimes{\rm id})\circ\Delta_{\mathcal{O}({\rm SU}_q(2))},
\qquad
i=a,c.
\end{gather*}
Using these formulas and the isomorphisms~\eqref{sphids}, one can check that the trivializable pieces of the comodule
algebra $\PWsq$ are isomorphic to $C(D)\otimes \mathcal{O}({\rm SU}_q(2))$ (cf.~\cite[equation~(1.8)]{hrz}).

Furthermore, as the comodule algebra $\PWsq$ is a~cotensor pro\-duct, combining~\cite[Lemma~1.13]{hrz}
with~\cite[Theorem~1.5]{hrz} yields that $\PW$ is a~piecewise trivial $(\operatorname{Ker}\pi)$-reduction
(see~\eqref{mudef}) of~$\PWsq$.

\begin{theorem}
[main result] The comodule algebra $\PWsq$ is not isomorphic to any smash product $C(S^2)\rtimes \mathcal{O}({\rm SU}_q(2))$
comodule algebra.
\end{theorem}

\pagebreak

\begin{proof}
Suppose that there exists a~cleaving map
\begin{gather*}
\mathcal{O}({\rm SU}_q(2)) \longrightarrow  \PWsq
\end{gather*}
that is an algebra homomorphism.
It is tantamount to the existence of a~${\rm U}(1)$-equivariant algebra homomorphism $f\colon \mathcal{O}({\rm SU}_q(2))\rightarrow\PW$
\cite[Proposition~4.1]{bz12}.
Let $\alpha$ and $\gamma$ denote generators of $\mathcal{O}({\rm SU}_q(2))$, and $a$, $c$ their classical counterparts.
Since $f([\alpha,\alpha^*])=0$, it follows from~\eqref{ac} that $f(\gamma)=0$ and $f(\alpha)f(\alpha)^*=1$.

On the other hand, by the ${\rm U}(1)$-equivariance, $ f(\alpha)=f_1a+f_2c$, for some $f_1,f_2\in C(S^2)$.
Furthermore, any continuous section of the Hopf line bundle $L_1$ can be written as $g_1a+g_2c$ for some $g_1,g_2\in
C(S^2)$.
We can rewrite it as $(g_1a+g_2c)f(\alpha)^*f(\alpha)$.
Since $(g_1a+g_2c)f(\alpha)^*\in C(S^2)$, we conclude that $f(\alpha)$ spans $\Gamma(L_1)$ as a~left $C(S^2)$-module.
Also, if $gf(\alpha)=0$ for some $g\in C(S^2)$, then $g=gf(\alpha)f(\alpha^*)=0$.
Hence $f(\alpha)$ is a~basis of $\Gamma(L_1)$ contradicting its nonfreeness.
\end{proof}

\subsection*{Acknowledgements} 

The authors are grateful to Tomasz Brzezi\'nski for discussions, and to the referees for
careful proofreading of the manuscript.
This work was partially supported by the NCN-grant 2011/01/B/ST1/06474.

\pdfbookmark[1]{References}{ref}
\LastPageEnding


\begin{thebibliography}{99}
\footnotesize\itemsep=0pt

\bibitem{bdh}
Baum P.F., De~Commer K., Hajac P.M., Free actions of compact quantum groups of
  unital $C^*$-algebras, \href{http://arxiv.org/abs/1304.2812}{arXiv:1304.2812}.

\bibitem{bh}
Baum P.F., Hajac P.M., Local proof of algebraic characterization of free
  actions, \href{http://arxiv.org/abs/1402.3024}{arXiv:1402.3024}.

\bibitem{bhms07}
Baum P.F., Hajac P.M., Matthes R., Szyma\'nski W., Noncommutative geometry
  approach to principal and associated bundles, in Quantum Symmetry in
  Noncommutative Geometry, {t}o appear, \href{http://arxiv.org/abs/math.DG/0701033}{math.DG/0701033}.

\bibitem{bh04}
Brzezi{\'n}ski T., Hajac P.M., The {C}hern--{G}alois character,
  \href{http://dx.doi.org/10.1016/j.crma.2003.11.009}{\textit{C.~R.~Math. Acad. Sci. Paris}} \textbf{338} (2004), 113--116,
  \href{http://arxiv.org/abs/math.KT/0306436}{math.KT/0306436}.

\bibitem{bz12}
Brzezi{\'n}ski T., Zieli{\'n}ski B., Quantum principal bundles over quantum
  real projective spaces, \href{http://dx.doi.org/10.1016/j.geomphys.2011.12.008}{\textit{J.~Geom. Phys.}} \textbf{62} (2012),
  1097--1107, \href{http://arxiv.org/abs/1105.5897}{arXiv:1105.5897}.

\bibitem{g-r99}
G{\"u}nther R., Crossed products for pointed {H}opf algebras, \href{http://dx.doi.org/10.1080/00927879908826704}{\textit{Comm.
  Algebra}} \textbf{27} (1999), 4389--4410.

\bibitem{hkmz11}
Hajac P.M., Kr{\"a}hmer U., Matthes R., Zieli{\'n}ski B., Piecewise principal
  comodule algebras, \href{http://dx.doi.org/10.4171/JNCG/88}{\textit{J.~Noncommut. Geom.}} \textbf{5} (2011), 591--614,
  \href{http://arxiv.org/abs/0707.1344}{arXiv:0707.1344}.

\bibitem{qsng}
Hajac P.M., Matthes R., So\l{}tan P.M., Szyma\'nski W., Zieli\'nski B.,
  Hopf--Galois extensions and $C^*$ algebras, in Quantum Symmetry in
  Noncommutative Geometry, {t}o appear.

\bibitem{hrz}
Hajac P.M., Rudnik J., Zieli\'nski B., Reductions of piecewise trivial comodule
  algebras, \href{http://arxiv.org/abs/1101.0201}{arXiv:1101.0201}.

\bibitem{kn63}
Kobayashi S., Nomizu K., Foundations of dif\/ferential geometry. {V}ol.~{I},
  Interscience Publishers, New York~-- London, 1963.

\bibitem{p-p95}
Podle{\'s} P., Symmetries of quantum spaces. {S}ubgroups and quotient spaces of
  quantum {${\rm SU}(2)$} and {${\rm SO}(3)$} groups, \href{http://dx.doi.org/10.1007/BF02099436}{\textit{Comm. Math.
  Phys.}} \textbf{170} (1995), 1--20, \href{http://arxiv.org/abs/hep-th/9402069}{hep-th/9402069}.

\bibitem{s-p99}
Schauenburg P., Galois objects over generalized {D}rinfeld doubles, with an
  application to {$u_q({\mathfrak{sl}}_2)$}, \href{http://dx.doi.org/10.1006/jabr.1998.7814}{\textit{J.~Algebra}} \textbf{217}
  (1999), 584--598.

\bibitem{s-pm11}
So{\l}tan P.M., On actions of compact quantum groups, \textit{Illinois~J.
  Math.} \textbf{55} (2011), 953--962, \href{http://arxiv.org/abs/1003.5526}{arXiv:1003.5526}.

\bibitem{w-sl87}
Woronowicz S.L., Twisted {${\rm SU}(2)$} group. {A}n example of a
  non-commutative dif\/ferential calculus, \href{http://dx.doi.org/10.2977/prims/1195176848}{\textit{Publ. Res. Inst. Math. Sci.}}
  \textbf{23} (1987), 117--181.

\end{thebibliography}
\end{document}